\documentstyle{article}
\def\m@th{\mathsurround=0pt }
\renewcommand{\Box}{\mbox{$\vrule height6pt width6pt depth0pt$}}
\def\eqalign#1{\null\,\vcenter{\openup\jot \m@th
   \ialign{\strut\hfil$\displaystyle{##}$&$
      \displaystyle{{}##}$\hfil \crcr#1\crcr}}\,}
\newtheorem{proposition}{Proposition}[section]
\newtheorem{theorem}{Theorem}[section]
\newtheorem{lemma}{Lemma}[section]
\newtheorem{corollary}{Corollary}[section]
\newtheorem{definition}{Definition}

\begin{document}
\begin{center}
{\Large{\bf A Ring-Theorist's Description of \\ Fedosov Quantization}}
\vskip 5mm
DANIEL R. FARKAS \\
Department of Mathematics \\
Virginia Polytechnic Institute and State University \\
Blacksburg, VA 24061
\end{center}
\vskip 1cm
\begin{quote}
{\footnotesize
ABSTRACT \ We present a formal, algebraic treatment of Fedosov's 
argument that the coordinate algebra of a symplectic manifold has a 
deformation quantization. His remarkable formulas are established in 
the context of affine symplectic algebras. 
}
\end{quote}

One of the standard techniques for reducing a problem about a noncommutative
algebra $B$ to one about commutative rings is to find a filtration $B_0 \subseteq
B_1 \subseteq B_2 \subseteq \cdots $ with $\bigcup B_i = B$ and $B_i B_j \subseteq
B_{i+j}$ so that the associated graded algebra is commutative. There is a shadow of
noncommutativity which can be recovered inside $gr B$: if $\overline{a}\in gr_m B$
and $\overline{b}\in gr_n B$, pull back $\overline{a}$ and $\overline{b}$ to $a$
and $b$ in $B_m$ and $B_n$ respectively. Then the image of the commutator $\lbrack a,b
\rbrack$ in
$gr_{m+n-1} (B)$ depends only on $\overline{a}$ and $\overline{b}$. This new binary
operation in $gr B$ is an example of a Poisson bracket.

In general, a (commutative) Poisson algebra $A$ over the scalar field $k$ is a
commutative $k$-algebra which is, at the same time, a Lie algebra under a second
operation $\lbrace \ast ,\ast \rbrace :A \times A \rightarrow A$ such that the
``product" or ``Liebniz" rule holds,
$$\lbrace ab,c \rbrace = a \lbrace b,c \rbrace + \lbrace a,c \rbrace b \quad {\rm for 
\ all} \  a,b,c \in A.$$
It is well known that the example $gr B$ described above is a Poisson algebra. It
can be generalized via the Rees construction to a special case of the following
artifact. Suppose the algebra $C$ has a central non-zero-divisor $h$ such that $C
/ hC$ is a commutative algebra.  Then $C / hC$ is a Poisson algebra under the
bracket $\lbrace \overline{c_1},\overline{c_2} \rbrace = {1 \over h}
\overline{\lbrack {c_1},{c_2} \rbrack }$. (Here the overbar denotes the canonical
homomorphism.)

A Poisson algebra has a ``naive quantization" if the process can be reversed; 
given a Poisson algebra $A$, there
is an algebra $C$ as we've just described such that $C / hC$ is isomorphic to $A$.
One can ask for much more.  The Poisson algebra $A$ has a
``deformation" or ``$\ast$-product" quantization provided there is a formal
deformation of the multiplication in $A$,
$$(a,b) \mapsto ab+ \lbrace a,b \rbrace h + \sum_{t=2}^{\infty} \mu_t (a,b)h^t$$
which satisfies the traditional associative identity and $\mu_t (a,b) = (-1)^t
\mu_t (b,a)$.

One of the earliest profound results about formal quantization of Poisson algebras is
the theorem of DeWilde and Lecomte ([DL]) asserting that if $M$ is a symplectic
manifold then $C^{\infty} (M)$ has a $\ast$-product quantization. In 
the 1990s, Fedosov
([Fe]) produced a new proof which he described as a ``simple geometric
construction". The purpose of the article which follows is to convince the reader
that Fedosov's construction is fundamentally algebraic. Whether it is actually
simple we leave to others to decide. Of course, the the theorem in 
question has been superseded rather spectacularly by Kontsevich's 
description of deformations for arbitrary Poisson manifolds. The 
author is not aware of any clear relationship between this new 
construction and Fedesov's.

Our context for interpreting Fedosov's techniques is the study of commutative
affine algebras over a field $k$ of characteristic zero, which are symplectic
algebras.  In earlier work, Loose ([L]), Farkas, and Letzter ([FL]) discussed
possible definitions of a symplectic algebra which either generalize or are
analogues of $C^{\infty} (M)$. In the ``smooth" case that $A$ is a regular affine
domain, all of the proposed approaches agree. Since we limit ourselves to these
algebas, we can define $A$ to be symplectic provided it is also a Poisson algebra
whose Poisson structure produces sufficiently many derivations. In detail, for
each $a \in A$ let $ham \, a$ denote the $k$-linear derivation of $A$ which sends $b$
to $\lbrace a,b \rbrace$. Then $A$ is {\it symplectic} when these derivations
generate all $k$-linear derivations as an $A$-module. (We use the shorthand $Der
\, A = A \cdot Ham \, A$.) Following Fedosov, we shall prove that every regular
affine symplectic domain has a naive quantization. (In fact, the noncommutative
algebra which is exhibited is a deformation quantization, which the interested
reader will be able to verify with ease.)

In very broad strokes, the naive quantization arises as the infinitesimal fixed
ring of a universal enveloping algebra. The Poisson bracket on $A$ induces an
$A$-bilinear alternating form $\omega$ on $Der \, A$, the module of $k$-derivations
from $A$ to itself, by extending $\omega (ham \, a,ham \, b)=\lbrace a,b \rbrace$.
Form the universal enveloping algebra $U$ with coefficients in $A$ based on the Lie
algebra $Ah \oplus Der \, A$ where $h$ is central and $\lbrack X,Y \rbrack = h
\omega (X,Y)$ for derivations $X$ and $Y$. This construction appears frequently in
the literature as an appropriate generalization of the Weyl algebra, cf. [B],[NR].
Roughly speaking, there is a derivation $D$ of $U$ for which $\lbrace x \in U |
D(x)=0 \rbrace$ turns out to be the quantization of $A$. It is, perhaps, worth
noting that nothing so elaborate is required if we restrict attention to graded
symplectic algebras; they all appear as the associated graded algebra for the ring
of differential operators on some regular domain ([F1]).

The mathematical reality is that the enveloping algebra is too small on two counts.
First, it must be completed and, secondly, we throw all differential $n$-forms
with coefficients into this larger algebra. The derivation $D$ can then be identified
with a connection as generalized by Fedosov. Most of our work is devoted to setting
up an algebraic apparatus for this ring so that some marvelous formulas of Fedosov
can be invoked. In particular, the Poisson bracket for a symplectic algebra induces a
duality between $Der \, A$ and the module of K\"ahler differentials on $A$. Fedosov
extends the duality to a ``homotopy inverse'' relation of differentials on the DGA
of $n$-forms just described. This can be done with ``local'' computations in the
geometric context, which are not available algebraically. Instead, we have to jump
back and forth between $U$ and the copy of its symmetrization which injects into $U$.
(The reader is not expected to know any symplectic geometry in order to follow this
paper. We hope the open-minded reader will
recognize the potential of algebra to clarify some mathematical ideas and to open
the door to abstractions which allow other problems to be resolved.) 

This paper was submitted in the early months of 1995 to a collection 
honoring the memory of S. A. Amitsur. Shimshon 
was especially encouraging in his support of my project to develop symplectic themes
algebraically. Unfortunately, inappropriate action of an editor led to 
the demise of the article. It was judged moot upon the final 
acceptance of [D]. I was surprised and delighted when Professor 
Sternheimer suggested in November 1999 that the manuscript be 
resurrected in its present form. This version contains the proofs 
referred to or abbreviated in the published paper of the same title.

\section[1]{The rings}

For the duration of this paper $k$ is a field of characteristic zero and $A$ is a
commutative affine domain which is regular and symplectic. The regularity assumption
will appear in the following form: the $A$-module of $k$-linear K\"ahler
differentials, $\Omega^1 (A)$, is a finitely generated projective $A$-module. We
will find it useful to exhibit projective bases explicitly. There is a standard
duality (cf. [MR]) between $\Omega^1 (A)$ and $Der \, A$, the module of $k$-linear
derivations from $A$ to itself. If $Z\in Der \, A$ and $u\,dv \in \Omega^1 (A)$ then
the pairing is given by $\langle Z,u \, dv \rangle = uZ(v)$. We may choose a
projective basis for $\Omega^1 (A)$ of the form
$$\lbrace (da_1 , \langle X_1 ,- \rangle ) , (da_2 , \langle X_2 ,- \rangle ) ,
\ldots ,(da_n ,\langle X_n , - \rangle ) \rbrace$$
where $a_1 , \ldots ,a_n \in A$ and $X_1 , \ldots ,X_n \in Der \, A$. In this case
we have what we have termed an {\it elite} basis for $Der \, A$,
\begin{equation}Y=\sum_j Y(a_j)X_j \quad {\rm for \ every} \ Y\in Der \, A.
\label {1} \end{equation}

The $A$-module $\Omega^q (A)$ consisting of (alternating) differential $q$-forms is
also finitely generated projective. We will write $\Omega (A)$ for the finite
direct sum $\sum_{q \geq 0} \Omega^q (A)$.

As discussed in [FL], we may define a regular affine domain $A$ to be symplectic
when it is a Poisson algebra with $Der \, A = A \cdot Ham \, A$. In other words,
$Der \, A$ is generated by $\lbrace ham \, a | \,  a \in A \rbrace$ as an
$A$-module. For any Poisson algebra there is an $A$-module map $\sharp : \Omega^1
(A) \rightarrow Der \, A$ given by $(\sum a_i \, db_i )^{\sharp} = \sum a_i \, ham
\, b_i$ ([H]). To say that $A$ is symplectic is to assert that $\sharp$ has an
inverse, which we denote $\flat$; indeed, a surjective map between finitely
generated projective modules of the same rank must be an isomorphism. Now define
$\omega : Der\, A \times Der\, A \rightarrow A$ by $\omega (X,Y) = \langle
X,Y^{\flat} \rangle$. Then $\omega$ is an $A$-bilinear form. It is nondegenerate,
for if $X\in Der\, A$ and $\omega (X,Y) = 0$ for all derivations $Y$ then $\langle
X,da\rangle = 0$ for all $a\in A$. But $X(a)=0$ for all $a$ forces $X=0$. Less
obviously, the form is alternating. It suffices to check this when $X=ham\, a$ and
$Y=ham\, b$.
$$\omega (ham\, a,ham\, b)=\langle ham\, a,db\rangle = \lbrace a,b\rbrace = -\lbrace
b,a\rbrace = -\omega (ham\, b,ham\, a).$$
We can now obtain a second formula for derivations, given an elite basis as
described above. If $Y\in Der\, A$ then
$$Y^{\flat}=\sum_j \langle X_j,Y^{\flat} \rangle da_j.$$ Hence
\begin{equation}Y=\sum_j \omega (X_j ,Y)ham\, a_j \quad {\rm for \ all} \ Y\in Der
\,A.
\label {2} \end{equation}

Let $h$ be an indeterminate, so that $hA$ represents a free module of rank one. Set
$L=hA\oplus Der\, A$, another projective $A$-module. More crucially, $L$ is a Lie
algebra over $A$ under the bracket
$$\lbrack (ha,X),(hb,Y)\rbrack = h\omega (X,Y)$$
for $a,b\in A$ and $X,Y\in Der\, A$. Equivalently, $h$ is central in $L$ and
$\lbrack X,Y\rbrack = h\omega (X,Y)$. We warn the reader of a treacherous difficulty
with this notation. The commutator bracket on $Der\, A$ which arises 
from composition is entirely different. For
the few times it will appear, we will use the notation $\lbrack X_{\circ} Y\rbrack$
for the derivation $X\circ Y-Y\circ X$. Returning to $L$, we can regard it as a
graded Lie algebra by assigning nonzero members of $Der\, A$ the degree 1 and to $h$ the
degree 2. Form the universal enveloping algebra $U_A (L)$. Since $L$ is a projective
$A$-module, $U_A (L)$ has a Poincar\'e-Birkhoff-Witt theorem in the sense that the
associated graded algebra for the filtration by products in $L$ is the symmetric
$A$-algebra $S_A (L)$ ([R]). As in the case of enveloping algebras over
fields ([D], 2.4), there is a homogeneous symmetrization map which injects $S_A
(L)$ into $U_A (L)$ by sending $y_1 y_2 \cdots y_m$ in $S_A (L)$ with $y_j \in L$ to
$${1\over m!} \sum_{\sigma \in {\cal S}_m} y_{\sigma (1)} y_{\sigma (2)} \cdots
y_{\sigma (m)}$$
in $U_A (L)$.

We are really interested in the subalgebra $S_A (Der\, A)$ of $S_A (L)$. This is
due to the fact that factoring out $h$ yields a surjective $A$-algebra homomorphism
$\pi : U_A (L) \rightarrow S_A (Der\, A)$. (The assertion is a consequence of the
universal properties of both $U_A (L)$ and $S_A (Der\, A)$.) Moreover, if $\pi
^\dagger :S_A (Der\, A) \rightarrow U_A (L)$ denotes the restriction of the
symmetrization map then $\pi ^\dagger$ is an $A$-module homomorphism with $\pi \circ
\pi ^\dagger = id$. We denote by $S_m ^\dagger$ the image in $U_A (L)$ of the ${\rm
m}^{\rm th}$ homogeneous component of $S_A (Der\, A)$ under $\pi ^\dagger$.

By applying induction on the length of words in $L$ to the formula $\lbrack
X,Y\rbrack = h\omega (X,Y)$ for $X,Y \in Der\, A$ we find that $\lbrack x,y \rbrack
\in hU_A (L)$ for all $x,y \in U_A (L)$. Consequently, ${1\over h} \lbrack x,y
\rbrack$ makes sense.  This leads to a weak quantization result.
\begin{theorem}The form $\omega$ extends uniquely to an $A$-bilinear Poisson bracket
$\{ *,*\}_{\omega}$ on $S_A (Der\, A)$ in such a way that for all $x,y \in U_A (L)$,
$$\pi ({1\over h} \lbrack x,y\rbrack )=\lbrace \pi (x),\pi (y)\rbrace _{\omega } .$$
\end{theorem}
\medskip
\noindent
{\bf Proof:}  It is easy to verify that the formula above can be used to define a
binary operation on $S_A (L)$. It is also obvious that the resulting bracket is
$A$-linear and satisfies the appropriate Liebniz and Lie requirements. Moreover, if
$X,Y\in Der\, A$ then
$$\lbrace X,Y\rbrace _{\omega} = \lbrace \pi (X),\pi (Y)\rbrace _{\omega} = \pi
({1\over h} \lbrack X,Y\rbrack ) = \pi \omega (X,Y) = \omega (X,Y).$$
Since $Der\, A$ generates $S_A (Der\, A)$ as an $A$-algebra, the bracket's
behavior on $Der\, A$ determines its behavior on the entire algebra.{\ \ \Box}  
\medskip

Since $L$ is a graded Lie algebra, $U_A (L)$ is a graded associative $A$-algebra.
This grading is compatible, under $\pi$ and $\pi ^\dagger$, with the natural grading
on $S_A (Der\, A)$. As a permanent notation we use $W$ for the completion of $U_A
(L)$ and $S$ for the completion of $S_A (Der\, A)$ with respect to these gradings.
With a slight abuse of notation we still have the completed maps $\pi$ and $\pi
^\dagger$. (These functions commute with infinite sums. Without further notice, the
term ``homomorphism'' will include this continuity property whenever it makes
sense.) Let $W_p$ denote the ${\rm p}^{\rm th}$ homogeneous component of $W$ and
write $W=\overline {\sum}_p W_p$ for the complete direct sum. Set $S^\dagger =
\overline {\sum}_p S_p ^\dagger$, an $A$-submodule of $W$. Observe that $S_p
^\dagger \subseteq W_p$.  More precisely, $\pi$ and $\pi ^\dagger$ preserve degree.

The quantization of $A$ turns out to be the $k$-subalgebra of $W$ on which a
certain derivation vanishes. This derivation, however, is defined on a slightly
larger ring. Consider $W\otimes _A \Omega (A)$ and its shadow $S\otimes _A \Omega
(A)$. (From now on we simply write $\Omega$ instead of $\Omega (A)$.) Both will be
regarded as $A$-algebras in the sense of ${\bf Z}_2$-graded algebras. Thus if $x\in
W\otimes _A \Omega ^{q(1)}$ and $y\in W\otimes _A \Omega ^{q(2)}$ then the graded
commutator is defined to be
$$\lbrack x,y \rbrack = xy - (-1) ^{q(1)q(2)} yx .$$ 
The center of $W\otimes \Omega$ consists of those elements $z$ such that $\lbrack
z,t \rbrack = 0$ for all $t\in W\otimes \Omega$. This definition is designed so that
$\Omega$ lies in the center of $W\otimes \Omega$. Another easy consequence of the
definition of graded bracket is that 
$$\lbrack x,y \rbrack \in h(W\otimes \Omega )$$
for all $x,y\in W\otimes \Omega$. Finally, recall that a graded derivation $D$ from
$W\otimes \Omega$ to itself of weight $N$ is a $k$-linear (continuous) map such that
$$D(W\otimes \Omega ^q ) \subseteq W\otimes \Omega ^{q+N}$$ for all $q$, and
$$D(xy) = D(x)y + (-1)^q xD(y)$$
for $x\in W\otimes \Omega ^q$ and arbitrary $y$.

Notice that $W\otimes \Omega$ is bigraded with homogeneous components $W_p \otimes
\Omega ^q$. We say that an element in this submodule has $W$-degree $p$ and
$\Omega$-degree $q$. Since $z$ lies in the center of $W\otimes \Omega$ if and only if
$\lbrack X,z\rbrack = 0$ for all $X\in Der\, A$, we see that each $W_p \otimes
\Omega ^q$ component of a central element is central.

Extend $\pi$ (resp. $\pi ^\dagger$) to a map from $W\otimes \Omega$ to $S\otimes
\Omega$ (resp. $S\otimes \Omega$ to $W\otimes \Omega$) by identifying it with $\pi
\otimes id$ (resp. $\pi ^\dagger \otimes id$).
\begin{theorem}Every element of $W\otimes \Omega$ can be written uniquely in the
form $\overline {\sum}_{m\geq 0} \,  h^m s_m$ with $s_m \in S^\dagger \otimes \Omega$.
\end{theorem}
\medskip
\noindent
{\bf Proof} We first treat existence. It suffices to show by induction on $p$ that
every element $u$ of $W_p$ has the form $\sum_{0\leq m <p} h^m s_m$ with $s_m \in
S_{p-2m} ^\dagger$. If $p=0$ then $W_0 = A$, so $u\in S_0 ^\dagger$; if $p=1$ then
$W_1 = Der\, A$, in which case $u\in S_1 ^\dagger$. In general,
$$u-\pi ^\dagger \pi (u) \in W_p .$$
Since $\pi (u-\pi ^\dagger \pi (u)) = 0$, we may write
$$u-\pi ^\dagger \pi (u) = hv$$
for some $v\in W_{p-2}$. Apply induction to $v$.

As to uniqueness, we must show that if $\overline {\sum} _{m\geq 0}\  h^m s_m = 0$
with $s_m \in S^\dagger \otimes \Omega$ then $s_m = 0$ for all $m$. Since $U_A (L)$
is an integral domain (e.g., from the fact that $A$ is a domain and the PBW
theorem), $W\otimes \Omega$ has no $h$-torsion. Thus if the sum is not formally
zero, we may assume that $s_0$ is not zero.
$$0 = \pi (0) = \pi \left( \overline{\sum} h^m s_m \right) = \pi (s_0 ).$$ 
But $\pi$ is injective on $S^\dagger$. Hence $s_0 = 0$, a contradiction. {\ \ \Box} 
\medskip
\begin{theorem}Assume that $\Phi : W\otimes \Omega \rightarrow W\otimes \Omega$ and
$\Phi _S : S\otimes \Omega \rightarrow S\otimes \Omega$ are graded derivations which
satisfy $$\pi \circ \Phi = \Phi _S \circ \pi.$$ If $\Phi (Der\, A) \subseteq \Omega$
and $\Phi (\Omega ) \subseteq \Omega$ then $$\Phi \circ \pi ^\dagger = \pi ^\dagger
\circ \Phi _S .$$
\end{theorem}
\medskip
\noindent
{\bf Proof:} We will repeatedly use the observation that the commutative square for
$\pi$ and $\Phi$ implies that $\Phi$ and $\Phi _S$ are formally the same on both
$Der\, A$ and $\Omega$. It suffices to prove the desired equality when restricted
to $S_p$ for all $p$. Indeed, if $c\in S_p$ and $\mu \in \Omega$ then
$$\eqalign{\left( \Phi \circ \pi ^\dagger \right) (c\otimes \mu ) & = \Phi \left( \pi
^\dagger (c)\mu \right) \cr
& = \left( \Phi \circ \pi ^\dagger \right) (c)\mu +(-1)^p \pi ^\dagger (c)\Phi
(\mu) \cr
& = \left( \pi ^\dagger \circ \Phi _S \right) (c)\mu +(-1)^p \pi ^\dagger \left(
c\, \Phi _S (\mu )\right) \cr
& = \left( \pi ^\dagger \circ \Phi _S \right) (c \otimes \mu) .}$$ 

So now consider $Y_1 ,\ldots ,Y_m \in Der\, A$ and fix $r\in \lbrace 1,\ldots
,m\rbrace$. For each choice $\sigma \in {\cal S}_m$ exactly one of the $m$ terms in
$\Phi (Y_{\sigma (1)} \cdots Y_{\sigma (m)} )$ has a $\Phi (Y_r)$ appearing. Thus if
we expand $\Phi \left( \sum _{\sigma \in {\cal S}_m} Y_{\sigma (1)} \cdots Y_{\sigma
(m)} \right)$ we will get $m!$ terms with $\Phi (Y_r)$ as a factor. Now $\Phi (Y_r)$
commutes with all derivations of $A$ so we can write any such term as $P\Phi (Y_r)$
with $P$ a product of the $Y_i \,$'s with $Y_r$ omitted. For each order of the factors
in $P$ we get $m$ repetitions of this term (i.e., one for each original placement of
$\Phi (Y_r)$ ). In other words, the contribution to
$$\Phi \left( {1\over m!} \sum_{\sigma} Y_{\sigma (1)} \cdots Y_{\sigma (m)}
\right)$$
of terms with the form $P\Phi (Y_r)$ is the symmetrization of $Y_1 \cdots \hat{Y_r}
\cdots Y_m$ multiplied by $\Phi (Y_r)$. The restricted formula follows.{\ \ \Box}
\medskip
\begin{theorem}The center of $W\otimes \Omega$ is $A[[h]]\otimes_A \Omega$.
\end{theorem}
\medskip
\noindent
{\bf Proof:} We already know that $A[[h]]\otimes_A \Omega$ is contained in the
center. So suppose $\xi$ is a central element. Write
$$\xi = \overline {\sum _{m\geq 0}} h^m \xi _m$$
with $\xi _m \in S^\dagger \otimes \Omega$. We first argue that each $\xi _m$ lies
in the center. This will follow if we show that $[Y,\xi _m]=0$ for all $Y\in Der\,
A$. The strategy is to use Theorem 1.1 : ${1\over h}ad\, Y$ is a graded derivation on
$W\otimes \Omega$ and $ham_{\omega}Y$ is a graded derivation on $S\otimes \Omega$
such that $$\pi \circ {1\over h}ad\, Y=ham_{\omega} Y\circ \pi.$$ 
Moreover, $\left( {1\over h}ad\, Y\right) \left( Der\, A\right) \subseteq \Omega ^0$
and $\left( {1\over h}ad\, Y\right) \left( \Omega \right) = 0$. Now apply Theorem
1.3.
$$\left( {1\over h}ad\, Y\right) \circ \pi ^\dagger =\pi ^\dagger \circ
ham_{\omega} Y.$$ In particular,
$$\left( {1\over h}ad\, Y\right) \left( S^\dagger \otimes \Omega \right) \subseteq
S^\dagger \otimes \Omega.$$ Thus
$${\overline {\sum_{m\geq 0}}} h^m \left( {1\over h}ad\, Y\right) \left( \xi _m
\right)$$
is the unique $h$-series expansion for $\left( {1\over h}ad\, Y\right) \left( \xi
\right)$. On the other hand, $\left( {1\over h}ad\, Y\right) \left( \xi 
\right) =0$. Therefore each $\left( {1\over h}ad\, Y\right) \left( \xi _m
\right) =0$. It follows that $\xi _m$ is central.

We have reduced the theorem to showing that if $\xi$ is central and lies in $S_p
^\dagger \otimes \Omega ^q$ then $\xi \in \Omega$. By applying induction to
formula (2) (with respect to an elite basis $(X_1 ,a_1 ), \ldots ,(X_n ,a_n )$ for
$Der\, A$), we may prove that if $Y_1, \dots ,Y_p \in Der\, A$ then  
$$\sum_{j=1} ^n \lbrack Y_1 Y_2 \cdots Y_p , X_j \rbrack ham\, a_j \equiv pY_1 Y_2
\cdots Y_p \ \ \ \  (mod\  h) .$$
Thus if $u \in S_p ^\dagger$ then
$$\sum_{j=1} ^n \lbrack u, X_j \rbrack ham\, a_j \equiv pu \ \ \ \ (mod\ h) .$$
(Note that this vacuously covers the case $p=0$.) Direct calculation shows that if
$\nu \in \Omega ^q$ we can extend the congruence,
$$\sum_{j=1} ^n \lbrack u\otimes \nu , X_j \rbrack ham\, a_j \equiv p(u\otimes \nu ) \
\ \ \ (mod\ h).$$
If we apply this to $\xi$ we obtain $$0\equiv p\xi \ \ \ \ (mod\ h) .$$
In other words, either $p=0$ or $\pi (\xi )=0$. Of course, $\pi$ is injective on
$S^\dagger \otimes \Omega$. We conclude that $p=0$.{\ \ \Box}

\section[2]{The differentials}

The foundation of Fedosov's calculations is the existence of two differentials on
$W\otimes \Omega$, extending $\sharp$ and $\flat$, which are
``homotopy inverses'' of each other. It will take us some technical preparation to
establish an algebraic version.

Suppose that $\alpha :A\rightarrow \Omega$ is a derivation into an $A$-bimodule
and that
$f:Der\, A \rightarrow W\otimes \Omega$ is a $k$-linear transformation such that
$$f(aX)=\alpha (a)X+af(X),\quad {\rm and}$$
$$h\alpha \left( \omega (X,Y) \right) = \lbrack f(X),Y\rbrack + \lbrack X,f(Y)
\rbrack $$
for all $X,Y\in Der\, A$ and $a\in A$. Then $f$ and $\alpha$ extend to a single
$k[h]$-linear derivation from $U_A (L)$ to $W\otimes \Omega$. The proof is a
straightforward modification of the upper triangular trick which can be found in [J],
page 154.

We first apply this to the situation in which $\alpha =0$ and $f=\flat$. If $X\in
Der\, A$ then $X^\flat \in \Omega ^1$. Thus $X^\flat$ commutes with all elements of
$W$ (and, in particular, with any $Y\in Der\, A$). The equations above are satisfied
trivially. Hence $\flat$ extends to a derivation
$$\delta : U_A (L) \rightarrow W\otimes \Omega . $$
Since $\delta (W_p)\subseteq W_{p-1} \otimes \Omega ^1$ it is easy to extend
$\delta$ to a derivation
$$\delta : W\rightarrow W\otimes \Omega .$$
Finally, it extends to a graded derivation of weight 1 on all of $W \otimes
\Omega$ by sending $u\otimes \nu$ to $\delta (u)\nu$. (The calculation rests on the
observation that $\delta (W)\subseteq W\otimes \Omega ^1$, so that if $\nu \in
\Omega ^q$ then $\nu \delta (u')=(-1)^q \delta (u')\nu$.) We call this derivation
$\delta$ as well. Its crucial properties are that
$$\eqalign{&\delta (\Omega )=0\quad ; \cr
&\delta (h)=0\quad ; \quad {\rm and}\cr
&\delta (X)=X^\flat \quad {\rm for\ all}\ X\in Der\,A .}$$
With respect to the grading, $\delta (W_p \otimes \Omega ^q )\subseteq W_{p-1}
\otimes \Omega ^{q+1}$.

In addition, $\delta ^2 =0$. Indeed, it certainly suffices to check that $\delta ^2
(W_p )=0$ for all $p$. Now if $u,v\in W$ then
$$\eqalign{\delta ^2 (uv)&=\delta \left( \delta (u)v+u\delta (v)\right) \cr
&=\left( \delta ^2 (u)v-\delta (u)\delta (v)\right)+\left( \delta (u)\delta
(v)+u\delta ^2 (v)\right) \cr}$$
because $\delta (u)\in W\otimes \Omega ^1$ and $\delta$ is a graded derivation. We
are reduced to checking that $\delta ^2 (Der\,A)=0$. But this is immediate.

There is an induced derivation $\delta _S : S\otimes \Omega \rightarrow S\otimes
\Omega$ which is $A$-linear and sends $X$ in $Der\,A$ to $X^\flat$. It is built so
that $\pi \circ \delta = \delta _S \circ \pi$. We also have $\delta \circ \pi
^\dagger = \pi ^\dagger \circ \delta _S$ by virtue of Theorem 1.3.

The second differential will be based on a function which is defined on $\Omega$ and
vanishes on $W$. We will need the classical contraction map on $\Omega (A)$. If
$X\in Der\,A$ then $\iota _X : \Omega ^q (A) \rightarrow \Omega ^{q-1} (A)$ is an
$A$-linear map given by
$$(\iota _X \xi )(Y_1, \ldots ,Y_{q-1} )=\xi (X,Y_1 ,\ldots ,Y_{q-1} ) .$$
The convention is that $\Omega ^{-1} =0$ and
$$\iota _X (\nu) = \langle X,\nu \rangle \quad {\rm for} \ \nu \in \Omega ^1.$$
Explicitly,$$ \iota _X(a_0 da_1 da_2 \cdots da_q ) = \sum _{t=1} ^q (-1)^{t-1} a_0
da_1 da_2 \cdots d\breve a_t \cdots da_q $$
where $d\breve a_t$ indicates that this factor is replaced with $X(a_t)$. It is
standard that $\iota_X$ is a graded derivation from $\Omega$ to itself of weight -1.
\begin{lemma}There is a unique $A$-linear map
$$\delta^* : \Omega (A) \rightarrow Der\,A \otimes _A \Omega (A)$$
such that $\delta ^* (\nu )=\nu ^\sharp $ for $\nu \in \Omega ^1$ and
$$\delta ^* (\alpha \beta ) = \delta ^* (\alpha )\beta +(-1)^q \alpha \delta ^*
(\beta )$$
for $\alpha \in \Omega ^q$ and arbitrary $\beta \in \Omega$.
\end{lemma}
\medskip
\noindent
{\bf Proof:} Recall that we have a pairing on $Der\,A \otimes \Omega ^1$ obtained
from the $A$-module isomorphism
$$Der\,A \simeq Hom_A (\Omega ^1 ,A).$$
Since $\Omega ^q$ is a projective $A$-module, we may tensor and obtain
$$Der\,A \otimes _A \Omega ^q \simeq Hom_A (\Omega ^1 ,\Omega ^q ).$$
This is isomorphism can be described by an extension of the original pairing: if
$X\in Der\,A,\  \nu \in \Omega ^q$, and $\lambda \in \Omega ^1$ define
$$\langle X\otimes \nu ,\lambda \rangle '=\langle X,\lambda \rangle \nu.$$
Now if $\phi \in Hom_A (\Omega ^1 ,\Omega ^q )$ there exists a unique $f\in Der\,A
\otimes \Omega ^q$ such that $\phi (\lambda )=\langle f,\lambda \rangle '$ for all
$\lambda \in \Omega ^1$. The new pairing has an additional attractive property. If
$\mu \in \Omega ^r$ and $g\in Der\,A\otimes \Omega ^q$ then
$$\mu \langle g,\lambda \rangle '=\langle \mu g,\lambda \rangle '\quad {\rm and}\quad
\langle g,\lambda \rangle '\mu = \langle g\mu ,\lambda \rangle '.$$

Fix $\nu \in \Omega ^{q+1}$. Consider the function which sends $\lambda$ in $\Omega
^1$ to $-\iota _{\lambda ^\sharp} (\nu)$; it lies in $Hom_A (\Omega ^1 ,\Omega
^q)$. Thus there exists an element we will call $\delta ^* (\nu )$ in $Der\,A
\otimes \Omega ^q$ such that
$$\langle \delta ^* (\nu ),\lambda \rangle '=-\iota _{\lambda ^\sharp} \nu \quad
{\rm for \ all}\ \lambda \in \Omega ^1 .$$
Clearly $\delta ^*$ is $A$-linear; we posit $\delta ^* (A)=0$. If $\nu \in \Omega
^1$ then
$$\eqalign{\langle \delta ^* (\nu ),\lambda \rangle &=-\iota _{\lambda ^\sharp} \nu
\cr
&=-\langle \lambda ^\sharp ,\nu \rangle \cr
&=-\omega (\lambda ^\sharp ,\nu ^\sharp ) \cr
&=\omega (\nu ^\sharp ,\lambda ^\sharp ) \cr
&=\langle \nu ^\sharp ,\lambda \rangle .}$$
Hence $\delta ^*$ agrees with $\sharp$ on $\Omega ^1$.

Finally, we shall prove that $\delta ^*$ has the derivation-like property because
$\iota _{\lambda ^\sharp}$ is a graded derivation. We borrow the symbols $\alpha$
and $\beta$ from the statement of the lemma and assume $\beta$ is homogeneous.
$$\eqalign{\langle \delta ^* (\alpha \beta ),\lambda \rangle '&=\iota _{\lambda
^\sharp} (\alpha \beta ) \cr
&=\iota _{\lambda ^\sharp} (\alpha )\beta +(-1)^q \alpha \iota _{\lambda ^\sharp }
(\beta ) \cr
&=\langle \delta ^* (\alpha ),\lambda \rangle '\beta +(-1)^q \alpha \langle \delta
^* (\beta ),\lambda \rangle ' \cr
&=\langle \delta ^* (\alpha )\beta +(-1)^q \alpha \delta ^* (\beta ),\lambda \rangle
'.{\ \ \Box}}$$
\medskip
\begin{lemma}$(\delta _S ^* )^2 =0$.
\end{lemma}
\medskip
\noindent
{\bf Proof:} It suffices to prove that $(\delta _S ^* )^2 (\Omega ^q )=0$ for all
$q$. We proceed by induction on $q$. The equality is true for $q=0$ and $q=1$ by
degree considerations. The induction step amounts to showing that if $\nu \in \Omega
^q$, $(\delta _S ^* )^2 (\nu )=0$ and $b\in A$ then $(\delta _S ^* )^2 (\nu\,
db)=0$. 
$$\eqalign{(\delta _S ^* )^2 (\nu\, db)&=\delta _S ^* (\delta _S ^* (\nu )db+(-1)^q
\nu
\delta _S ^* (db)) \cr
&=(\delta _S ^* )^2 (\nu )db+(-1)^{q-1} \delta _S ^* (\nu )\delta _S ^* (db)+(-1)^q
\delta _S ^* (\nu )\delta _S ^* (db)+\cr
&\quad \quad \quad \quad \quad \quad (-1)^{2q} (\delta _S ^* )^2 (db) \cr
&=0. {\ \ \Box}}$$
\medskip
\begin{lemma}If $x\in S_p \otimes \Omega ^q$ then
$$(\delta _S \delta _S ^* +\delta _S ^* \delta _S )(x)=(p+q)x.$$
\end{lemma}
\medskip
\noindent
{\bf Proof:} We first check this for $x\in S_p$. Since $\delta _S ^*$ vanishes on
$S$, we require $\delta _S ^* \delta _S (Y_1 \cdots Y_p )=pY_1 \cdots Y_p$ for $Y_j
\in Der\,A$. Calculate:
$$\eqalign{\delta _S ^* \delta _S (Y_1 \cdots Y_p )&=\delta _S ^* \left( \sum
_{j=1} ^p Y_1 \cdots \hat Y_j \cdots Y_p Y_j ^\flat \right) \cr
&= \sum
_{j=1} ^p Y_1 \cdots \hat Y_j \cdots Y_p Y_j ^{\flat \sharp} \cr
&=pY_1 \cdots Y_p .}$$
(Though the role of commutativity is crucial here, we will see in the next lemma how
to lift the calculation to $W$.)

Next we plug in $\mu \in S^q$. This time we recall that $\delta _S (\mu)=0$. The
formula $\delta _S \delta _S ^* (\mu) =qu$ is proved inductively on $q$. For $q=1$,
$$\delta _S \delta _S ^* (\mu )=\delta _S (\mu ^\sharp )=\mu ^{\sharp \flat} =\mu .$$
If $\delta _S \delta _S ^* (\mu) =qu$ and $b\in A$ then
$$\eqalign{\delta _S \delta _S ^* (\mu \,db)&=\delta _S \left( \delta _S ^* (\mu
)db+(-1)^q \mu (db)^\sharp \right) \cr
&=\delta _S \delta _S ^* (\mu )db+(-1)^{2q} \mu (db)^{\sharp \flat} \cr
&=(q+1)\mu \,db .}$$

Last of all, we evaluate $\delta _S \delta _S ^* +\delta _S ^* \delta _S$ on
$y\otimes \nu$ for $y\in S_p$ and $\nu \in \Omega ^q$.
$$\eqalign{&\quad \left( \delta _S \delta _S ^* +\delta _S ^* \delta _S \right)
(y\otimes \nu ) \cr
&=\delta _S \left( y\delta _S ^* (\nu )\right) + \delta _S ^* \left( \delta _S
(y)\nu \right) \cr
&=\delta _S (y)\delta _S ^* (\nu )+y\delta _S \delta _S ^* (\nu ) +\delta _S ^*
\delta _S (y)\nu +(-1)\delta _S (y)\delta _S ^* (\nu ) \cr
&=qy\nu + py\nu .{\ \ \Box}}$$
\medskip

We shall jack Lemma 2.3 up to $W\otimes \Omega$. It turns out that the formula in
the lemma characterizes $S^\dagger$.
\begin{proposition}Let $x\in W_p \otimes \Omega ^q$. Then $x\in S_p ^\dagger \otimes
\Omega ^q$ if and only if
$$(\delta ^* \delta + \delta \delta ^* )(x)=(p+q)x .$$
\end{proposition}
\medskip
\noindent
{\bf Proof:} To prove the formula for $x\in S_p ^\dagger \otimes \Omega ^q$ it
suffices to verify it when $x\in S_p ^\dagger $ --- the rest of the proof of Lemma
2.3 holds for $\delta ^*$ and $\delta$ replacing $\delta _S ^*$ and $\delta _S$
respectively. So we may assume that $x={1\over p!}\sum _{\sigma \in {\cal S}_p}
Y_{\sigma (1)} \cdots Y_{\sigma (p)}$ with $Y_1 ,\ldots ,Y_p$ in $Der\,A$. For each
$\sigma \in {\cal S}_p$
$$\eqalign{\delta ^* \delta (Y_{\sigma (1)} \cdots Y_{\sigma (p)} )&=\delta ^*
\left( \sum_{j=1} ^p Y_{\sigma (1)} \cdots \hat Y_{\sigma (j)} \cdots Y_{\sigma
(p)} Y_{\sigma (j)} ^\flat \right) \cr
&=\sum_{j=1} ^p Y_{\sigma (1)} \cdots \hat Y_{\sigma (j)} \cdots Y_{\sigma
(p)} Y_{\sigma (j)} .}$$ 
Thus each original permutation of $Y_1 ,\ldots ,Y_p$ gives rise to $p$ permutations,
each ending with a different choice of $Y_t$. From the opposite point of view, if we
take any product ending in $Y_t$, the last letter could have moved to the end
position (via the $\delta ^*$ shuffle) from any one of $p$ positions. Thus each
permutation appears $p$ times as we list the terms in the expansion of $\delta ^*
\delta \left( \sum _{\sigma \in {\cal S}_p}
Y_{\sigma (1)} \cdots Y_{\sigma (p)} \right)$. In other words, $\delta ^* \delta
(x)=px$.

Conversely, assume that $x\in W_p \otimes \Omega ^q$ satisfies the formula. Write
$x=\overline{\sum}_{m\geq 0} \,h^m s_{p-2m}$ where $s_t \in S_t ^\dagger \otimes
\Omega ^q$. Then
$$(\delta ^* \delta + \delta \delta ^* )(x)=\overline{\sum _{m\geq 0}} h^m (p+q-2m)
s_{p-2m}$$ 
according to the previous direction of the proposition. By the uniqueness of the
$h$-series,
$$p+q=p+q-2m$$ for all $m$ with $s_{p-2m} \neq 0$. That is, $m=0$.{\ \ \Box}
\medskip

We conclude this section with Fedosov's fundamental ``homotopy inverse'' relation.
Unfortunately, this will involve one further layer of notation because $\delta ^*$
does not quite meet all of his requirements. We perturb $\delta ^*$ in defining a
new function $\tilde \delta : W\otimes \Omega \rightarrow W\otimes \Omega$.

If $a\in A$ set $\tilde \delta (a)=0$. If $x\in S_p ^\dagger \otimes \Omega ^q$ for
$(p,q)\neq 0$ set
$$\tilde \delta (x)={1\over p+q} \pi ^\dagger \delta_S ^* \pi .$$
Extend $\tilde \delta \ A[[h]]$-linearly to all of $W\otimes \Omega$ using the
$h$-series. Note that $\tilde \delta (W_p \otimes \Omega ^q ) \subseteq W_{p+1}
\otimes \Omega ^{q-1}$ and $\tilde \delta (W)=0$. We still have $\tilde \delta (\nu
)=\nu ^\sharp$ for $\nu \in \Omega ^1$.

Consider the augmentation map from $S$ to $A$. It has two relevant extensions to all
of $S\otimes \Omega$. The first, denoted $\tau _S : S\otimes \Omega \rightarrow A$,
is the composition of the projection $S\otimes \Omega \rightarrow S$ with the
augmentation map. The second, denoted $T_S : S\otimes \Omega \rightarrow \Omega$,
is the augmentation map tensored with the identity. Thus $T_S = (\tau _S |S)\otimes
id$. Both $T_S$ and $\tau_S$ are $A$-linear maps which are, in a sense, idempotent. 

If $x\in W\otimes \Omega$, decompose $x=\overline{\sum}_{m\geq 0} \, h^m s_m$ with
$s_m \in S^\dagger \otimes \Omega$ and define $\tau : W\otimes \Omega \rightarrow
A[[h]]$ by 
$$\tau (x)=\overline {\sum_{m\geq 0}} h^m \tau _S \pi (s_m) .$$
Here $\tau _S \pi$ can also be regarded as the composition of the projection
$W\otimes \Omega \rightarrow W$ with the augmentation map $W\rightarrow A$. The
uniqueness of the $h$-series implies that $\tau$ is $A[[h]]$-linear. Set $T :
W\otimes \Omega \rightarrow A[[h]]\otimes _A \Omega$ to be $(\tau |W)\otimes id$. We
remark that $\tau$ is somewhat mysterious. For example, if $X$ and $Y$ are
derivations of $A$ then $XY={1\over 2}(XY+YX)+{1\over 2}[X,Y]$. Hence
$$\tau (XY)={1\over 2}h\omega (X,Y) .$$

It immediately follows from the definitions that
$$\pi \circ \tau =\tau _S \circ \pi \quad {\rm and}\quad \pi \circ T=T_S \circ \pi
.$$
\begin{theorem}
\begin{enumerate}
\item[({\it i})]
$\tilde \delta ^2=0$.
\item[({\it ii})]
(Fundamental Formula)$\quad \tilde \delta \delta +\delta \tilde \delta +\tau = id$.
\end{enumerate}
\end{theorem}
\medskip
\noindent
{\bf Proof:} Since all of the maps in the statement of the theorem are
$A[[h]]$-linear, it suffices to verify the equalities when restricted to $S_p
^\dagger \otimes \Omega ^q$. Let $x$ be in this component.
$$\eqalign{\tilde \delta \left( \tilde \delta (x)\right) &= \tilde \delta \left(
{1\over p+q}
\pi ^\dagger \delta _S ^* \pi (x)\right) \cr
&=\left( {1\over p+q} \right)^2 \pi ^\dagger \delta _S ^* \pi \pi ^\dagger \delta _S
^*
\pi (x) .}$$
Since $\pi \pi ^\dagger = id$ and $(\delta _S ^* )^2 =0$ (by Lemma 2.2), the last
expression is zero. Thus $\tilde \delta ^2 (x)=0$.

As to ({\it ii}), if $(p,q)=0$ then $x\in A$. In this case, $\delta (x)=\tilde
\delta (x)=0$. Since $\tau |A =id$,
$$\tilde \delta \delta (x) +\delta \tilde \delta (x) + \tau (x)=x.$$
If $(p,q)\neq 0$ then $\tau (x)=0$.
$$\eqalign{(\tilde \delta \delta +\delta \tilde \delta )(x)&={1\over p+q} \left( \pi
^\dagger \delta _S ^* \pi \delta +\delta \pi ^\dagger \delta _S ^* \pi \right) (x)
\cr
&={1\over p+q} \left( \pi
^\dagger \delta _S ^* \delta _S \pi +\pi ^\dagger \delta _S \delta _S ^* \pi \right)
(x) }$$
by the commuting squares for $\delta$ with $\pi$ and $\pi ^\dagger$. Invoking Lemma
2.3 we obtain
$$(\tilde \delta \delta +\delta \tilde \delta )(x)=\pi ^\dagger \pi (x) .$$
Since $x\in S^\dagger\otimes \Omega$ we have $\pi ^\dagger \pi (x)=x$.{\ \ \Box}

\section[3]{Connections}

A connection for $A$ associates to each pair of derivations $X,Y\in Der\,A$ a third
derivation $\nabla _X Y$ with the properties that $\nabla$ is $A$-linear in the
lower argument and
\begin{equation}\nabla _X \,aY=X(a)Y+a\nabla _X Y
\label {3} \end{equation}
for all $a\in A$. We will initially think of $\nabla Y$ as a function from $Der\,A$
to itself via $(\nabla Y)(X)=\nabla _X Y$. Given an elite basis $(X_1 ,a_1 ),\ldots
,(X_n ,a_n)$ for $Der\,A$ then
$$\nabla _X Y=\nabla_{\sum X(a_i )X_i} Y=\sum X(a_i)\nabla _{X_i} Y .$$
In this way we can identify $\nabla Y$ with $\sum \nabla _{X_i} Y\otimes da_i$ in
$Der\,A\otimes \Omega ^1$ and regard $\nabla$ as a function $Der\,A \rightarrow
Der\,A \otimes \Omega ^1$.

It is tempting to try to extend $\nabla$ to a $k$-endomorphism of $W\otimes \Omega$.
The classical insight is to define $\nabla a$ to be $da$ for $a\in A$. Then (3)
becomes
$$\nabla aY=(\nabla a)Y+a\nabla Y .$$
We are now on familiar ground. The discussion at the beginning of the second section
tells us that $\nabla$ extends to a $k[[h]]$-linear derivation from $W$ into
$W\otimes \Omega ^1$ provided
$$h\,d\omega (X,Y) = [\nabla X,Y]+[X,\nabla Y]$$
in $W\otimes \Omega ^1$.

Recall that a connection $\nabla$ for $A$ is {\it parallel} (to the symplectic form
$\omega$) when
$$Z\left( \omega (X,Y)\right) =\omega (\nabla _Z X,Y)+\omega (X,\nabla_Z Y)$$
for all $X,Y,Z \in Der\,A$. It is a {\it Poisson} {\it connection} if it also
satisfies the ``no torsion'' condition
$$\nabla _X Y-\nabla _Y X=[X_{\circ} Y] .$$
An extensive algebraic discussion of Poisson connections can be found in [F2]. In
particular, a regular symplectic affine domain such as $A$ always supports a Poisson
connection.
\begin{lemma}If $\nabla$ is a parallel connection then
$$h\,d\omega (X,Y)=[\nabla X,Y]+[X,\nabla Y]$$
for all $X,Y \in Der\,A$.
\end{lemma}
\medskip
\noindent
{\bf Proof:} We expand $[\nabla X,Y]+[X,\nabla Y]$ using the formula
$$\nabla V=\sum \left( \nabla _{X_i} V\right) da_i$$
for $V=X$ and $V=Y$.
$$\eqalign{[\nabla X,Y]+[X,\nabla Y]&=\sum _i h\left( \omega (\nabla _{X_i}
X,Y)+\omega (X,\nabla _{X_i} Y)\right) da_i \cr
&=\sum _i hX_i (\omega (X,Y)) da_i \cr
&=h\sum _i \langle X_i ,d\omega (X,Y)\rangle da_i \cr
&=h\, d(\omega (X,Y)) .{\ \ \Box}}$$
\medskip

Assume from now on that $\nabla$ is a Poisson connection. The induced derivation
$\nabla : W\rightarrow W\otimes \Omega ^1$ extends to a graded derivation $\nabla :
W\otimes \rightarrow W\otimes \Omega$ of weight 1 by setting
$$\nabla (x\otimes \nu )=(\nabla x)\nu +x\,d\nu .$$
(That the extension is well-defined follows from the centrality of $a$ and $da$ for
$a\in A$.) The Poisson connection has a similar extension $_S \nabla : S\otimes
\Omega \rightarrow S\otimes \Omega$ which is related to $\nabla$ by the commutative
square $_S \nabla \circ \pi =\pi \circ \nabla$. We will need several less trivial
identities.
\begin{lemma}If $\nabla$ is a Poisson connection for $A$ then
$$\omega (X,\nabla _Z ham\,u)=\omega (Z,\nabla _X ham\,u)$$
for all $X,Z\in Der\,A$ and $u\in A$.
\end{lemma}
\medskip
\noindent
{\bf Proof:}$$\omega (X,\nabla _Z ham\,u)=Z(X(u))-\omega (\nabla _Z X,ham\,u)$$
because $\nabla$ is parallel. Since $\nabla$ has no torsion,
$$\eqalign{\omega (\nabla _Z X, ham\,u)&=\omega (\nabla _X Z,ham\,u)+\omega
([Z_\circ X],ham\,u) \cr
&=\omega (\nabla _X Z,ham\,u)+Z(X(u))-X(Z(u)) .}$$
Hence
$$\eqalign{\omega (X,\nabla _Z ham\,u)&=X(Z(u))-\omega (\nabla _X Z,ham\,u) \cr
&=\omega (Z,\nabla _X ham\,u) }$$
by one more application of parallelism.{\ \ \Box}
\medskip
\begin{lemma}$\delta \nabla$ vanishes on $Ham\,A$.
\end{lemma}
\medskip
\noindent
{\bf Proof:} We re-examine the previous lemma with the help of a projective basis
$(da_1 ,X_1 ),\ldots ,(da_n ,X_n )$ for $\Omega ^1$. First of all, if $Y,Z\in
Der\,A$ and $u\in A$ then
$$\eqalign{\langle Y,\delta (\nabla _Z ham\,u )\rangle &=\langle Y,(\nabla _Z ham\,u
)^\flat \rangle \cr
&=\omega (Y,\nabla _Z ham\,u) \cr
&=\omega (Z,\nabla _Y ham\,u) }$$
with the last equality from Lemma 3.2. Now
$$\nabla _Y ham\,u=\sum _i Y(a_i )\nabla _{X_i} ham\,u .$$
Therefore
$$\langle Y,\delta (\nabla _Z ham\,u )\rangle =\langle Y,\sum _i \omega (Z,\nabla
_{X_i} ham\,u)da_i \rangle.$$
Consequently
$$\delta (\nabla _Z ham\,u)=\sum _i \langle Z,\delta (\nabla _{X_i}
ham\,u)\rangle da_i. \eqno(*)$$

Using the induced elite basis a second time, we recall that $\nabla ham\,u=\sum _j
(\nabla _{X_j} ham\,u)da_j$, so
$$\eqalign{\delta \nabla ham\,u&=\sum_j \delta (\nabla _{X_j} ham\,u)da_j \cr
&=\sum_j \sum_i \langle X_j ,\delta (\nabla _{X_i} ham\,u)\rangle da_i da_j \quad
{\rm by}\ (*) \cr
&=-\sum _i \left( \sum _j \langle X_j ,\delta (\nabla _{X_i} ham\,u)\rangle da_j
\right) da_i \cr
&=-\sum _i \delta (\nabla _{X_i} ham\,u)da_i \quad {\rm by\ the\ basis\ property}
\cr
&=-\delta \nabla ham\,u .{\ \ \Box}}$$
\medskip
\begin{theorem}$\delta \nabla +\nabla \delta =0$
\end{theorem}
\medskip
\noindent
{\bf Proof:} Notice that $(\delta \nabla +\nabla \delta )(\xi)=0$ for $\xi \in
\Omega$ because $\delta$ vanishes on $\Omega$. If $\xi =ham\,u$ for $u\in A$ then
$$(\nabla \delta )(ham\,u)=\nabla (du)=0 .$$
Thus Lemma 3.3 implies that $(\delta \nabla +\nabla \delta)(\xi)=0$ for this choice
of $\xi$ as well. The value of these observations is that $\Omega$ and $Ham\,A$
generate $W\otimes \Omega$ as a $k[[h]]$-algebra (in the sense that they
algebraically generate its homogeneous components).

The theorem will be completed once we prove that if $r$ and $t$ are bihomogeneous
elements of $W\otimes \Omega$ with
$$(\delta \nabla +\nabla \delta )(r)=(\delta \nabla +\nabla \delta )(t)=0$$
then
$$(\delta \nabla +\nabla \delta )(rt)=0 .$$
But a direct calculation shows that
$$(\delta \nabla +\nabla \delta )(rt)=(\nabla \delta (r)+\delta \nabla
(r))t+r(\nabla \delta (t)+\delta \nabla (t)) .{\ \ \Box}$$
\medskip
\begin{theorem}$T\circ \nabla =\nabla \circ T$
\end{theorem}
\medskip
\noindent
{\bf Proof:} By modifying the argument in Theorem 1.3 or Proposition 2.1 the reader
may verify that $\pi ^\dagger \circ \nabla ={}_S\nabla \circ \pi ^\dagger$. (The
bookkeeping will look familiar by writing $\nabla Y_t =\sum _j \nabla _{X_j} Y_t
\otimes da_j$ with respect to an elite basis.)

Since $T$ and $\nabla$ are both $k[[h]]$-linear, it suffices to prove that $T\circ
\nabla \circ \pi ^\dagger =\nabla \circ T\circ \pi ^\dagger$. A quick glance at the
definition shows that $T\circ \pi ^\dagger =T_S$. We are reduced to checking
that $T_S \circ {}_S \nabla = {}_S \nabla \circ T$. This is an immediate consequence
of the facts that
$$\eqalign{{}_S\nabla (S_p \otimes \Omega ^q )&\subseteq S_p \otimes \Omega ^{q+1}
\quad {\rm and} \cr 
 T_S | S_p \otimes \Omega ^q &=\left\{ \begin{array}{ll} 0 & \mbox{if $p\neq
0$}
\\
id & \mbox{if $p=0$ .}
\end{array} \right.  {\ \ \Box}}$$

\medskip
The last general property we will need about connections is really an observation
about curvature. It turns out that we do not rely on any geometric properties of
curvature other than that it exists in a very weak algebraic manner of speaking:
$\nabla ^2$ is essentially inner. This is a special case of a more general
phenomenon, based on the simple calculation that $\nabla ^2$ is a graded
$k[[h]]$-linear derivation of $W\otimes \Omega$ which vanishes on $\Omega$.
\begin{theorem}Let $D: W\otimes \Omega \rightarrow W\otimes \Omega$ be a graded
derivation such that $D(W_1 )\subseteq W_1 \otimes \Omega$. Then $D$ vanishes on
$A[[h]]\otimes \Omega$ if and only if there exists a $\Gamma \in W\otimes \Omega$
such that $D={1\over h}ad\,\Gamma$ (i.e., $D(\alpha )={1\over h}[\Gamma ,\alpha ]$
for all $\alpha \in W\otimes \Omega$ ).
\end{theorem}  
\medskip
\noindent 
{\bf Proof:} Since $A[[h]]\otimes \Omega$ is the center of $W\otimes \Omega$, any
inner derivation of $W\otimes \Omega$ vanishes on $A[[h]]\otimes \Omega$.
Conversely, assume that $D$ vanishes on the center and set
$$\Gamma =-{1\over 2}\sum _i (ham\,a_i )DX_i$$
with respect to an elite basis. It suffices to show that $D(Y)={1\over h}[\Gamma
,Y]$ for $Y\in Der\,A$. By definition,
$$[\Gamma ,Y]=-{1\over 2}\sum _i [ham\,a_i ,Y]DX_i -{1\over 2}\sum _i (ham\,a_i
)[DX_i ,Y] .$$
We examine the first term.
$$\eqalign{\sum_i [ham\,a_i ,Y]DX_i &=h\sum_i \omega (ham\,a_i ,Y)DX_i \cr
&=-h\sum _i Y(a_i )DX_i \cr
&=-hD\left( \sum_i Y(a_i )X_i \right) \quad {\rm because}\ D(A)=0, \cr
&=-hDY .}$$ 
Thus $$[\Gamma ,Y]={1\over 2} hDY-{1\over 2}\sum _i (ham\,a_i )[DX_i ,Y] .$$
We now simplify the second term.
$$\eqalign{\sum _i (ham\,a_i )[DX_i ,Y]&=\sum_i (ham\,a_i )D\left( [X_i ,Y]\right)
-\sum _i (ham\,a_i )[X_i ,DY] \cr
&=-\sum _i (ham\,a_i )[X_i ,DY]}$$
because $D(A)=0$. According to the hypothesis of the theorem, $DY$ is a sum of
expressions of the form $Z\otimes \xi$ with $Z\in Der\,A$ and $\xi\in \Omega$.
$$\eqalign{\sum _i (ham\,a_i )[X_i ,Z\otimes \xi ]&=\sum _i (ham\,a_i )[X_i
,Z]\otimes \xi \cr
&=hZ\otimes \xi \quad {\rm by}\ (2).}$$
Hence $$\sum _i (ham\,a_i )[DX_i ,Y]=-hDY .$$
We conclude that
$$[\Gamma ,Y]={1\over 2}hDY+{1\over 2}hDY=hDY .{\ \ \Box}$$
\medskip
\begin{corollary}Let $\nabla$ be a Poisson connection for $A$. Then there exists
$R\in W_2 \otimes \Omega ^2$ with $\nabla ^2 ={1\over h}ad\, R.${\ \ \Box}
\end{corollary}

\section[4]{The Fedosov calculus}

One of Fedosov's insights is to expand the notion of connection to any graded
derivation on $W\otimes \Omega$ of weight one which formally satisfies equation (3).
This creates a source of connections large enough to always include one with no
``curvature''.

We begin with the technical engine which drives Fedosov's calculus.
\begin{definition}A function $\Phi : W\otimes \Omega \rightarrow W\otimes \Omega$ is
{\bf monotone} provided
$$\Phi \left( \overline{\sum _{p\geq t}} W_p \otimes \Omega \right) \subseteq
\overline{\sum _{p\geq t}} W_p \otimes \Omega$$
for all $t\in {\bf N}$.
\end{definition}
\medskip
\begin{theorem}[Vanishing Theorem] Assume $\Phi$ is monotone. If $\beta \in W\otimes
\Omega$ satisfies
\begin{enumerate}
\item[(i)] $\delta (\beta )=\Phi (\beta )$,
\item[(ii)] $\tilde \delta (\beta )=0$, and
\item[(iii)] $\tau (\beta )=0$
\end{enumerate}
then $\beta =0$.
\end{theorem}
\medskip
\noindent
{\bf Proof:} Assume $\beta \neq 0$. Write $\beta =\overline{\sum}_{p\geq t} \, b_p$
with $b_p \in W_p \otimes \Omega$ and $b_p \neq 0$. By the fundamental formula of
Theorem 2.1,
$$\beta =\tilde \delta \delta (\beta )+\delta \tilde \delta (\beta )+\tau (\beta
)=\tilde \delta \Phi (\beta ) .$$
Hence $$\overline{\sum _{p\geq t}} b_p \in \overline{\sum _{p\geq t}} \tilde \delta
(W_p \otimes \Omega ) .$$
But $\tilde \delta (W_p \otimes \Omega )\subseteq W_{p+1} \otimes \Omega$. That is,
$$\overline{\sum _{p\geq t}} b_p \in \overline{\sum _{p\geq t}} W_{p+1} \otimes
\Omega .$$
We have reached the contradiction $b_t =0$. {\ \ \Box}
\medskip
\begin{definition}A {\bf Fedosov} {\bf connection} $D$ for $A$ is a graded
$k[[h]]$-linear derivation of weight one on $W\otimes \Omega$ such that
$$\eqalign{&D(a)=da \mbox{ for all } a\in A ; \cr
&D+\delta \mbox{ is monotone; and} \cr
&D^2 =0 .}$$
\end{definition}
\medskip
\begin{theorem}Assume $D$ is a Fedosov connection for $A$. Then for each $u\in
A[[h]]$ there exists a unique $b\in W$ such that $D(b)=0$ and $\tau (b)=u$.
\end{theorem}
\medskip
\noindent
{\bf Proof:} We first treat uniqueness. Suppose $D(b)=0$ and $\tau (b)=0$. Then
$\delta (b)=(D+\delta )(b)$ and $\tilde \delta (b) =0$ because $\tilde \delta$ is
zero on $W$. The Vanishing Theorem yields $b=0$. 

As to existence, define a sequence
$b_0 ,b_1 ,\ldots $ in $W$ by 
$$b_0 =u \quad {\rm and} \quad b_{m+1} =\tilde \delta (D+\delta )(b_m) .$$
Since $\tilde \delta (D+\delta )\left( \overline{\sum} _{p\geq t}\, W_p\otimes \Omega
\right) \subseteq \overline{\sum} _{p\geq t+1} W_p\otimes \Omega$, the series $\sum
_{j=0} ^\infty b_j$ converges. Call this sum $b$. By construction,
$$b=u+\tilde \delta (D+\delta )(b) .$$
Observe that
$$\eqalign{\tau (b)&=b-\tilde \delta \delta (b) \cr
&=u+\tilde \delta (D+\delta )(b)-\tilde \delta \delta \tilde \delta (D+\delta )(b)
\cr
&=u+\tilde \delta (D+\delta )(b)-\tilde \delta (id -\tilde \delta \delta -\tau
)(D+\delta )(b) \cr
&=u+\tilde \delta (D+\delta )(b)-\tilde \delta (D+\delta )(b)+\tilde \delta \tau
(D+\delta )(b) .}$$
But $(D+\delta )(b) \in W\otimes \Omega ^1$ and $\tau$ is zero on $W\otimes \Omega
^1$. Hence
$$\tau (b)=u .$$
We conclude that
$$b=\tau (b)+\tilde \delta (D+\delta )(b) .\eqno(*)$$
The argument is completed by showing that $D(b)=0$ via ($*$) and the Vanishing
Theorem.
$$\delta D(b)=(D+\delta )D(b)-D^2 (b)=(D+\delta )D(b) .$$
This establishes condition ({\it i}) of Theorem 4.1 for $D(b)$. By ($*$), $b=\tau
(b)+\tilde \delta \delta (b)+\tilde \delta D(b)$. With the fundamental formula, this
yields
$$\delta \tilde \delta (b)=\tilde \delta D(b) .$$
Since $\tilde \delta (W)=0$, we obtain ({\it ii}). Finally, $\tau D(b)=0$ because
$D(b)\in W\otimes \Omega ^1$.{\ \ \Box}
\medskip
\begin{corollary}Given a Fedosov connection $D$ for $A$, set
$$B=\{b\in W | D(b)=0 \},$$
a subalgebra of $W$ which contains $h$. Then the restriction of the augmentation map
$\varepsilon : W\rightarrow A$ induces an isomorphism
$$B/hB \simeq A$$ and ${1\over h}[b,b']$ is sent to $\{\varepsilon (b),\varepsilon
(b')\}$. (In other words, $B$ is a naive quantization of $A$.)
\end{corollary}
\medskip
\noindent
{\bf Proof:} We first argue that $hB=B\cap {\rm Ker}\varepsilon$. Since $\varepsilon
(h)=0$, one inclusion is obvious. So suppose $b\in B\cap{\rm Ker} \varepsilon$.
Using the $h$-series expansion,
$$b=\overline {\sum _{m\geq 0}} h^m s_m \quad {\rm with}\quad s_m \in S^\dagger ,$$
we see that $0=\varepsilon (b)=\varepsilon (s_0 )$. Hence
$$\tau (b)=\overline {\sum _{m\geq 1}} h^m \varepsilon (s_m )=h\overline {\sum
_{m\geq 1}} h^{m-1} \varepsilon (s_m ) .$$
Since $\overline{\sum}_{m\geq 1} \,h^{m-1} \varepsilon (s_m ) \in A[[h]]$, the
theorem tells us that there exists $c\in B$ with $\tau (c)=\overline{\sum}_{m\geq 1}
\,h^{m-1} \varepsilon (s_m )$. Therefore
$$hc\in B\quad {\rm and}\quad \tau (hc)=\tau (b).$$
The uniqueness assertion of the theorem establishes $hc=b$. This proves $b\in hB$.

The restriction of $\varepsilon$ to $B$ is injective. It is surjective by the
existence portion of the theorem: for each $u\in A$ there is an element
$u+\overline{\sum}_{i>0} \,b_i$ in $B$ with $\varepsilon \left( \overline{\sum}_{i>0}
\,b_i \right) =0$ by virtue of the $W$-degrees involved. We now have 
$$B/hB \simeq A .$$

For the Poisson bracket calculation, we may assume that $\tau (b)=\varepsilon (b)$
and $\tau (b')=\varepsilon (b')$. The $W_1$-component of $b$ is
$$\tilde \delta \nabla \varepsilon (b)=\tilde \delta (d\varepsilon (b))=ham\,
(\varepsilon (b)) .$$
It follows that the $W_2$-component of $[b,b']$ is
$$h\omega (ham\,\varepsilon (b),ham\,\varepsilon (b')) .$$
But $\omega (ham\,u,ham\,v)=\{u,v\}$ in general. Hence the $W_0$-component of
${1\over h}[b,b']$ is $\{ \varepsilon (b),\varepsilon (b')\}$. {\ \ \Box}

\section[5]{Existence of Fedosov connections}

The last step in [Fe] is the construction of a Fedosov connection. Here we follow
the original quite closely, replacing classical identities about curvature with
simple algebraic calculations. The overall strategy is transparent: given a Poisson
connection, find a correction term $\gamma \in \overline{\sum} _{p\geq 3} W_p
\otimes \Omega ^1$ so that $\nabla -\delta +{1\over h} ad\,\gamma$ is a Fedosov
connection. Notice that this expression is a sum of graded derivations of weight 1,
which makes it a derivation of the same type. Moreover,
$$(\nabla -\delta +{1\over h}ad\,\gamma )(a)=da \quad \mbox{for all}\ a\in A$$
because $\delta$ vanishes on $A$ and $A$ is central. Also, $\nabla +{1\over
h}ad\,\gamma$ is monotone what with $\nabla$ preserving $W$-degree and ${1\over h}
ad\,\gamma$ pushing up the lowest $W$-degree in the support of an element by at
least 1. This leaves the crucial issue of showing that the corrected graded
derivation has square zero.
\begin{lemma}For any $\gamma \in W\otimes \Omega ^1$,
$$(\nabla -\delta +{1\over h}ad\,\gamma )^2 ={1\over h}ad\left( -\delta (\gamma )
+R+\nabla (\gamma )+ {1\over h}\gamma ^2 \right)$$
(Here $\nabla ^2 ={1\over h}ad\,R$ as in Corollary 3.1 .)
\end{lemma}
\medskip
\noindent
{\bf Proof:}$$\eqalign{(\nabla -\delta +{1\over h}ad\,\gamma )^2 =\nabla ^2 &+\delta
^2 +{1\over h^2} (ad\,\gamma )^2 \cr
&-\delta \nabla -\nabla \delta \cr
&-{1\over h}\left( \delta (ad\, \gamma )+(ad\,\gamma )\delta \right) \cr
&+{1\over h}\left( \nabla (ad\,\gamma )+(ad\,\gamma )\nabla \right) .}$$
We analyze the right-hand side of the equation line by line. First, $\nabla ^2
={1\over h}ad\,R$ and $\delta ^2 =0$. Since $\gamma$ has $\Omega$-degree 1,
$(ad\,\gamma )^2 =ad(\gamma ^2)$. It is easy to check that $\gamma ^2 \in
h(W\otimes \Omega ^2 )$ so we may rewrite the first line as
$${1\over h}ad(R+{1\over h}\gamma ^2 ) .$$
The second line on the right-hand side is zero by Theorem 3.1. Direct calculations
yield
$$\eqalign{\delta (ad\,\gamma )+(ad\,\gamma )\delta &=ad(\delta (\gamma )) \quad {\rm
and} \cr
\nabla (ad\,\gamma )+(ad\,\gamma )\nabla &=ad(\nabla (\gamma )).}$$
The lemma follows. {\ \ \Box}
\medskip

The goal is to construct $\gamma \in \overline{\sum}_{p\geq 3} W_p \otimes \Omega
^1$ so that $-\delta (\gamma )+R+\nabla (\delta )+{1\over h}\gamma ^2$ is central.
This requires some knowledge about $R$ which we delineate below.
\begin{lemma}
\begin{enumerate}
\item[(i)]$\delta R=0$.
\item[(ii)]$\nabla R=dT(R)$, where $d$ is extended to a $k[[h]]$-linear map on
$A[[h]]\otimes \Omega$.
\end{enumerate}
\end{lemma}  
\medskip
\noindent
{\bf Proof:} ({\it i}) By definition, $\nabla ^2 X={1\over h}[R,X]$ for all $X\in
Der\,A$. Since $\delta \nabla +\nabla \delta =0$ we know that $\delta \nabla ^2
X=\nabla ^2 \delta X$. Now $\delta X\in \Omega ^1$, so $\nabla ^2 \delta X \in
d^2(\Omega ^1)$. Consequently, $\nabla ^2 \delta X=0$. We infer that
$$\delta ([R,X])=0 \quad {\rm for}\ X\in Der\,A .$$
On the other hand, $\delta [R,X]=[\delta R,X]+[R,\delta X]$. Again, $\delta X$ is
central, whence
$$[\delta R,X]=0 \quad \mbox{for all}\ X\in Der\,A .$$
It follows that $\delta R$ lies in the center of $W\otimes \Omega$.

But $R\in W_2 \otimes \Omega ^2$ implies $\delta R\in W_1 \otimes \Omega ^3$. We
have already identified the center of $W\otimes \Omega$ as $A[[h]]\otimes \Omega$.
Since $h$ has $W$-degree 2, there are no nonzero elements of the center with odd
$W$-degree. Therefore $\delta R=0$.

({\it ii}) This identity is a variation on the associative law.
$$\nabla ^2 (\nabla X)={1\over h}[R,\nabla X] \quad {\rm for}\ X\in Der\,A .$$
On the other hand, $\nabla ^2 (\nabla X)=\nabla ^3 X=\nabla (\nabla ^2 X)$, and
$$\nabla (\nabla ^2 X)=\nabla ({1\over h}[R,X])={1\over h}[\nabla R,X]+{1\over
h}[R,\nabla X] .$$
It follows that $[\nabla R,X]=0$, i.e., $\nabla R$ is central. Put in a different
form,
$$\nabla R=T(\nabla R) .$$
Applying Theorem 3.2, we obtain $\nabla R=\nabla (T(R))$. Finally, $\nabla$ agrees
with $d$ on $A[[h]]\otimes \Omega$. {\ \ \Box}
\medskip
\begin{theorem}There exists a $\gamma \in \overline{\sum} _{p\geq 3} W_p \otimes
\Omega ^1$ such that
$$\nabla -\delta +{1\over h}ad\,\gamma$$
has square zero.
\end{theorem}
\medskip
\noindent
{\bf Proof:} Inductively define $\gamma _m \in W_m \otimes \Omega ^1$ by
$$\gamma _3 = \tilde \delta R \quad {\rm and}$$
$$\gamma _t ={1\over h} \sum _{p+q-1=t} \tilde \delta (\gamma _p \gamma _q ) +\tilde
\delta \nabla \gamma _{t-1} \quad {\rm for}\ t>3.$$
Then $\gamma =\overline{\sum}_{t\geq 3} \,\gamma _t$ satisfies
$$\gamma =\tilde \delta R+{1\over h}\tilde \delta (\gamma ^2)+\tilde \delta \nabla
\gamma .$$
Since $\tilde \delta ^2 =0$, we trivially have
$$\tilde \delta \gamma =0 .\eqno(*)$$
Set $\beta =\delta \gamma +T(R)-R-\nabla \gamma -{1\over h}\gamma ^2$.

Since $\delta$ vanishes in $A[[h]]\otimes\Omega$, $\delta ^2 =0$, and $\delta R=0$,
$$\eqalign{\delta \beta &=-\delta \nabla \gamma -{1\over h}\delta (\gamma ^2 ) \cr
&=-\delta \nabla \gamma +{1\over h}[\gamma ,\delta \gamma ] \quad \mbox{because
$\gamma$ has $\Omega$-degree 1} \cr
&=\nabla \delta \gamma +{1\over h}[\gamma ,\delta \gamma ] \quad \mbox{by Theorem
3.1} \cr
&=\nabla \delta \gamma +{1\over h}[\gamma ,\beta -T(R)+R+\nabla \gamma +{1\over
h}\gamma ^2 ] \cr
&=\nabla \delta \gamma +{1\over h}[\gamma,\beta] +{1\over h}[\gamma,R] +{1\over
h}[\gamma,\nabla \gamma ] .}$$
We compare this formula to one for $\nabla \beta$:
$$\eqalign{\nabla \beta &=\nabla \delta \gamma +dT(R)-\nabla R-\nabla ^2 \gamma
-{1\over h}\nabla (\gamma ^2 ) \cr
&=\nabla \delta \gamma -{1\over h}[R,\gamma ]-{1\over h} \nabla (\gamma ^2 ) .}$$
But $\nabla \gamma^2 =(\nabla \gamma )\gamma -\gamma (\nabla \gamma )=[\nabla \gamma
,\gamma ]$. Hence
$$\delta \beta =\nabla \beta +{1\over h}[\gamma ,\beta ]. \eqno({\it i}) $$
In addition,
$$\eqalign{\tilde \delta \beta &=\tilde \delta \delta \gamma -\tilde \delta
(R+\nabla \gamma +{1\over h}\gamma ^2 ) \cr
&=\tilde \delta \gamma -\gamma  \quad \mbox{according to the defining relation for
$\gamma$}\cr 
&=\delta \tilde \delta \gamma -\tau (\gamma ) \quad \mbox{by the
fundamental formula} \cr
&=0-0 \quad \mbox{by ($*$)} .}$$
Summarizing,
$$\tilde \delta \beta =0 . \eqno({\it ii})$$

We have two-thirds of the hypotheses for the Vanishing Theorem. Obviously $\tau
(\beta )=0$ because $\beta$ has $\Omega$-degree 2. Finally, $\nabla +{1\over
h}ad\,\gamma$ is monotone. We conclude that $\beta =0$. In other words,
$$-\delta \gamma +R+\nabla \gamma +{1\over h}\gamma ^2$$
is the central element $T(R)$. {\ \ \Box}
\medskip
\begin{corollary}Fedosov connections exist. {\ \ \Box}
\end{corollary}
\medskip

When we put together Corollaries 4.1 and 5.1 we arrive at the conclusion that any
regular affine domain which is symplectic has a naive quantization.

\bigskip
\bigskip
\centerline{{\bf REFERENCES}}
\bigskip

\medskip
\noindent
[B] J.C. Baez, R-Commutative geometry and quantization of Poisson algebras, {\it
Advances in Math.}, {\bf 95} (1992), 61-91.
 
\medskip
\noindent
[DL] M. De Wilde and P.B.A. Lecomte, Formal deformations of the Poisson algebra
of a symplectic manifold and star-products. Existence, equivalence, derivations,
In: {\it Deformation Theory of Algebras and Structures and Applications} (eds. 
M. Hazewinkel and M. Gerstenhaber), NATO ASI Ser. C {\bf 247}, Kluwer Acad. 
Publ., (1988), 897-960.

\medskip
\noindent
[D] J. Dixmier, {\it Enveloping Algebras}, North-Holland, Amsterdam (1977).

\medskip
\noindent
[D] J. Donin, On the quantization of Poisson brackets, {\it Advances 
in Math.}, {\bf 127} (1997), 73-93.

\medskip
\noindent
[F1] D.R. Farkas, Characterizations of Poisson algebras, {\it Comm. in 
Algebra}, {\bf 23(12)} (1995), 4669-4686.

\medskip
\noindent
[F2] D.R. Farkas, Connections, Poisson brackets, and affine 
algebras, {\it Journal of Algebra}, {\bf 180} (1996), 757-777.

\medskip
\noindent
[FL] D.R. Farkas and G. Letzter, Ring theory from symplectic 
geometry, {\it Journal of Pure and Applied Alg.}, {\bf 125} (1998), 
155-190.    

\medskip
\noindent
[Fe] B.V. Fedosov, A simple geometrical construction of deformation quantization,
{\it Journal of Differential Geom.}, {\bf 40} (1994), 213-238.

\medskip
\noindent
[H] J. Huebschmann, Poisson cohomology and and quantization, {\it J. Reine 
Angew. Math.}, {\bf 408} (1990), 57-113.

\medskip
\noindent
[J] N. Jacobson, {\it Lie Algebras}, Interscience, New York, (1962).

\medskip
\noindent 
[L] F. Loose, Symplectic algebras and Poisson algebras, {\it Comm. in 
Algebra}, {\bf 21(7)} (1993), 2395-2416.

\medskip
\noindent
[MR] J.C. McConnell and J.C. Robson, {\it Noncommutative Noetherian Rings}, 
John Wiley and Sons, Chichester (1987).

\medskip
\noindent
[NR] Y. Nouaz\'e and P. Revoy, Sur les alg\`ebres de Weyl 
g\'en\'eralis\'ees, {\it Bull. Sci. Math.}, {\bf (2) 96} (1972), 27-47.

\medskip
\noindent
[R] G. Rinehart, Differential forms for general commutative algebras, {\it Trans.
Amer. Math. Soc.}, {\bf 108} (1963), 195-222.

\end{document}